\documentclass{amsart}
\usepackage[utf8]{inputenc}
\usepackage[T1]{fontenc}
\usepackage{amsmath,amsthm,amsfonts,amssymb}
\usepackage{microtype}
\microtypesetup{expansion=false}
\usepackage[hidelinks]{hyperref}
\usepackage{tikz}
\newtheorem{theorem}{Theorem}[section]
\newtheorem{corollary}[theorem]{Corollary}
\newtheorem{proposition}[theorem]{Proposition}
\newtheorem{lemma}[theorem]{Lemma}
\numberwithin{equation}{section}
\newcommand{\dd}{\,d}
\newcommand{\C}{\mathbb C}
\newcommand{\R}{\mathbb R}

\begin{document}
\title[Ratio and zero limits from recurrence coefficients]
{Elementary ratio and zero limits for varying recurrences}
\author{K. Castillo}
\address{CMUC, Department of Mathematics, University of Coimbra,
3000-143 Coimbra, Portugal}
\email{math@keniercastillo.com}
\author{S. Yakubovich}
\address{Department of Mathematics, Faculty of Sciences,
University of Porto, Campo Alegre st., 687, 4169-007 Porto, Portugal}
\email{syakubov@fc.up.pt}
\date{October 12, 2026}
\begin{abstract}
Two theorems of Loureiro and Van Assche on ratio and zero limits
provided part of the motivation for this note. Their elementary
nature permits direct proofs; Appendix~\ref{sec:power-known} gives
complete proofs of both for the canonical initial conditions,
using earlier continued-fraction identities.
We prove ratio limits and weak limits of normalized zero measures for
varying higher-order three-term recurrences with arbitrary coefficient
periods and independent limiting functions for their phases.
Joint local coefficient limits give ratio convergence without uniform
coefficient bounds. For zero measures, convergence in measure to finite
measurable profiles suffices, without boundedness or integrability;
the limiting measures may have noncompact support. Coefficients may
vanish, and real diagonal terms are allowed for ordinary three-term
recurrences. The proofs use interlacing, partial fractions with
nonnegative residues, local computations of moments, and stability
of zero counts under changes of recurrence coefficients. Truncation
removes the coefficient bounds from the zero theorem, whose proof is
independent of ratio asymptotics. Uniform bounds give the version with
compact support. Power-law limits and their explicit transform formula
follow as a special case.
\end{abstract}
\maketitle

\section{General ratio and zero limits}

Theorems~\ref{thm:ratios} and~\ref{thm:zeros} describe how limits of
varying recurrence coefficients determine polynomial ratios and zero
measures. The coefficient period $p$ is arbitrary, and each phase has
its own limiting function. Coefficients may vanish. Ratio convergence
requires joint local coefficient limits, without uniform bounds.
For zero measures, convergence in measure to finite measurable profiles
suffices, without boundedness or integrability.
Proposition~\ref{prop:compact-zeros} gives compact support under uniform
bounds. Proposition~\ref{prop:local-ratios} also gives a ratio limit
under local coefficient limits at a diverging index, without requiring
a limit for that index divided by $N$.

Our starting point was the precise conditions on blocks of ordinary
three-term recurrence coefficients for exact polynomial mappings
in the work of de Jesus and Petronilho
\cite[Theorem~2.1]{DeJesusPetronilho2010}, also used in the study
of zeros of sieved ultraspherical polynomials
\cite{CastilloDeJesusPetronilho2020}. Independent perturbations
need not preserve these conditions; we therefore asked which ratio
and zero limits can be obtained directly from the recurrence
coefficients. During the authors' discussions, Yakubovich recalled
the symmetric Appell families in \cite{LoureiroYakubovich2013}
and the preprint \cite{LoureiroVanAssche2026}.
These Appell families satisfy higher-order three-term recurrences
with positive coefficients.

The common power-law regime studied by Loureiro and Van Assche
\cite[Theorems~1.1--1.2]{LoureiroVanAssche2026} is a special case of
these results, with period $r$. Corollary~\ref{cor:original} recovers
both limits for canonical initial conditions, including the explicit
transform formula. Appendix~\ref{sec:power-known} gives a particularly
elementary deduction of this positive power-law case from earlier
continued-fraction identities.

For positive coefficients and zero diagonal terms, Ben Romdhane
\cite[Theorem~2.2 and its proof]{BenRomdhane2008} proves the
location and interlacing of the zeros needed for the probability
representations of the ratios below.
Delvaux and L\'opez \cite{DelvauxLopez2013} give further interlacing
results and treat arbitrary coefficient periods.
Aptekarev, Kalyagin, and Saff \cite{AptekarevKalyaginSaff2009}
study the asymptotics of higher-order three-term recurrences.
For third-order recurrences with a single continuous limiting function,
Aptekarev and Novokshenov
\cite[Theorems~1--2]{AptekarevNovokshenov2025} study ratio asymptotics
and zero measures. Their ratio theorem assumes uniform containment
of the zeros on a finite star.

Kuijlaars and Rom\'an \cite{KuijlaarsRoman2010} obtain averaged
zero limits for higher-term recurrences with real interlacing zeros.
After the substitution $u=z^{r+1}$, their theorem also covers the
positive cases with zero diagonal terms and $p\mid(r+1)$, under
\eqref{eq:bound} and the joint coefficient limits \eqref{eq:local},
when the limiting functions are positive for $t>0$ and extend
continuously to $t=0$, and the joint coefficient limits also hold
as $m,N\to\infty$ with $m/N\to0$.
For ordinary three-term recurrences, limits with continuous functions
are studied in \cite{KuijlaarsVanAssche1999}.
Varying asymptotically periodic coefficients are treated in
\cite{VanAssche1999Periodic}.
Kuijlaars and Serra Capizzano \cite[Theorem~2.1]{KuijlaarsSerra2001}
allow measurable limiting functions and convergence in measure,
without boundedness or integrability assumptions.
Theorem~\ref{thm:zeros} gives the corresponding conclusion for zero
measures normalized to have mass one when $r=p=1$, and extends it to
symmetric recurrences
of higher order, with arbitrary periods and independent limiting
functions for the residue classes.
The proofs use interlacing and partial fractions together with a finite
calculation: each fixed coefficient in the expansions of the ratios and
their logarithms depends on finitely many nearby recurrence coefficients.
Zero moments are averages of these polynomials. Related formal recursions
and moment expansions occur in
\cite{PetreolleSokalZhu2023,LopezGarciaProkhorov2024}.
For zero measures, Lemma~\ref{lem:counting} controls the effect of changing
recurrence coefficients, regardless of the magnitude of those changes.
Truncation then removes the coefficient bounds and proves
Theorem~\ref{thm:zeros} independently of ratio asymptotics.

Fix integers $r,p\ge1$, and put $d=r+1$ and
$\omega=e^{2\pi i/d}$. The recurrence has order $d$;
the number of coefficient phases $p$ need not equal $r$.
For each $N\ge1$, define monic polynomials by
\begin{equation}\label{eq:recurrence}
 q_{n+1,N}(z)=(z-b_{n,N})q_{n,N}(z)-c_{n-r,N}q_{n-r,N}(z),
 \quad n\ge0,
\end{equation}
where $q_{0,N}=1$, $q_{j,N}=0$ for $j<0$, and $c_{m,N}=0$ for $m<0$.
Throughout, $c_{m,N}\ge0$ and $b_{m,N}\in\R$. We require
$b_{m,N}=0$ when $r\ge2$; when $r=1$, arbitrary real diagonal
coefficients are allowed. Thus the latter case need not be symmetric.
For the compact zero limits we use the bound
\begin{equation}\label{eq:bound}
 C_T:=\sup_{N\ge1,\ 0\le m\le TN}
       \max\{|b_{m,N}|,c_{m,N}\}<\infty,\quad T>0.
\end{equation}
Write
$$
 \Sigma=\bigcup_{j=0}^{d-1}\omega^j[0,\infty),
 \quad \Sigma_B=\bigcup_{j=0}^{d-1}\omega^j[0,B].
$$
For $r=1$, these sets are $\R$ and $[-B,B]$.
Each phase corresponds to a residue class modulo $p$.
All phase indices below are taken modulo $p$.

\begin{theorem}[Ratio limits]\label{thm:ratios}
Under \eqref{eq:recurrence} and the stated sign conditions, suppose there are functions
$a_k:(0,\infty)\to[0,\infty)$ and
$h_k:(0,\infty)\to\R$, $0\le k<p$, such that
\begin{equation}\label{eq:local}
 \begin{aligned}
 &c_{m,N}\to a_k(t),\quad b_{m,N}\to h_k(t)\\[7pt]
 &\text{whenever }N\to\infty,\ m/N\to t>0,\ m\equiv k\pmod p.
 \end{aligned}
\end{equation}
Condition~\eqref{eq:local} already implies continuity of these functions.
Here $h_k=0$ if $r\ge2$. As $N\to\infty$, whenever $n/N\to t>0$ and
$n\equiv k\pmod p$,
\begin{equation}\label{eq:ratios}
 \begin{aligned}
 \frac{q_{n,N}(z)}{q_{n+1,N}(z)}&\to\varphi_k(z,t),\\[7pt]
 \frac{q_{n,N}(z)}{q_{n+p,N}(z)}&\to
 \Phi(z,t)=\prod_{k=0}^{p-1}\varphi_k(z,t).
 \end{aligned}
\end{equation}
locally uniformly on $\C\setminus\Sigma$; the second limit
does not require a fixed phase. The limits are nonzero and satisfy
\begin{equation}\label{eq:system}
 \frac1{\varphi_k(z,t)}
 =z-h_k(t)-a_{k-r}(t)\prod_{j=1}^r\varphi_{k-j}(z,t),
 \quad z\varphi_k(z,t)\to1\quad(z\to\infty).
\end{equation}
This normalization selects a unique analytic branch near infinity.
Both convergences are uniform when $n/N$ ranges over a compact
subinterval of $(0,\infty)$, with the limits evaluated at $t=n/N$.
The same uniform convergence holds after any fixed number of derivatives in $z$.
\end{theorem}

For vectors $a\in[0,\infty)^p$ and $h\in\R^p$, with $h=0$
if $r\ge2$, let $\mu_{a,h}$ denote the limiting normalized zero
measure of the recurrence with periodic coefficients
$c_m=a_{m\bmod p}$ and $b_m=h_{m\bmod p}$.
Its existence is proved below. Write $\varphi_k(z;a,h)$ for the
functions in \eqref{eq:system} with $a_k(t)=a_k$ and $h_k(t)=h_k$, and set
$\Phi(z;a,h)=\prod_{k=0}^{p-1}\varphi_k(z;a,h)$.
This product is analytic and nonzero on $\C\setminus\Sigma$,
with $z^p\Phi(z;a,h)\to1$ at infinity. Set
$$
 A_{k,N}(t)=c_{p\lfloor Nt/p\rfloor+k,N},
 \quad H_{k,N}(t)=b_{p\lfloor Nt/p\rfloor+k,N}.
$$

\begin{theorem}[Limits of zero measures]\label{thm:zeros}
Under \eqref{eq:recurrence} and the stated sign conditions, let
$a_k:(0,\infty)\to[0,\infty)$ and $h_k:(0,\infty)\to\R$,
$0\le k<p$, be measurable and finite almost everywhere,
with $h_k=0$ if $r\ge2$. Suppose that, for every $T>0$ and every $k$,
\begin{equation}\label{eq:measure}
 A_{k,N}\longrightarrow a_k,\quad H_{k,N}\longrightarrow h_k
 \quad\text{in measure on }(0,T).
\end{equation}
As $N\to\infty$, whenever $n/N\to T>0$, the normalized zero
measures, counting multiplicity, converge weakly:
\begin{equation}\label{eq:zero-average}
 \mu_{n,N}:=\frac1n\sum_{q_{n,N}(\zeta)=0}\delta_\zeta
 \quad\longrightarrow\quad
 \mu_T=\frac1T\int_0^T\mu_{a(t),h(t)}\dd t.
\end{equation}
On $\C\setminus\Sigma$,
\begin{equation}\label{eq:G}
 \begin{aligned}
 \frac1n\frac{q'_{n,N}(z)}{q_{n,N}(z)}&\to
 G_T(z):=\int\frac{\dd\mu_T(\zeta)}{z-\zeta}\\[7pt]
 &=-\frac1{pT}\int_0^T
       \frac{\partial_z\Phi(z;a(t),h(t))}{\Phi(z;a(t),h(t))}\dd t
 \end{aligned}
\end{equation}
locally uniformly. No coefficient bounds, integrability of $a_k,h_k$,
or limits at $t=0$ are required. The limiting measure need not have
compact support.
\end{theorem}

The proof is given in Appendix~\ref{sec:unbounded}.
We first prove the following version with bounded coefficients.

\begin{proposition}[Zero measures with compact support]\label{prop:compact-zeros}
Under \eqref{eq:recurrence}--\eqref{eq:bound}, let
$a_k:(0,\infty)\to[0,\infty)$ and $h_k:(0,\infty)\to\R$,
$0\le k<p$, be measurable, with $h_k=0$ if $r\ge2$. Suppose
\begin{equation}\label{eq:L1}
 \|A_{k,N}-a_k\|_{L^1(0,T)}+
 \|H_{k,N}-h_k\|_{L^1(0,T)}\to0,\quad T>0,\quad 0\le k<p.
\end{equation}
Then \eqref{eq:zero-average} holds whenever $n/N\to T>0$.
The measures $\mu_{n,N}$ for all sufficiently large $N$ and $\mu_T$
are supported on $\Sigma_B$, where one may take
$B=B_r(C_L)$ for any $L>T$, with
\begin{equation}\label{eq:root-bound}
 B_r(C)=\begin{cases}
 C+2\sqrt C,&r=1,\\[7pt]
 \dfrac{r+1}{r}(rC)^{1/(r+1)},&r\ge2.
 \end{cases}
\end{equation}
The convergence in \eqref{eq:G} is locally uniform on
$\C\setminus\Sigma_B$, using the analytic continuation of
$\Phi(z;a(t),h(t))$ to this set.
Under \eqref{eq:bound}, conditions~\eqref{eq:measure} and
\eqref{eq:L1} are equivalent.
\end{proposition}

Under \eqref{eq:bound}, condition~\eqref{eq:local} gives pointwise
convergence of $A_{k,N}$ and $H_{k,N}$; dominated convergence then
gives \eqref{eq:L1}.
The next corollary recovers both main theorems of
\cite{LoureiroVanAssche2026} by taking $p=r$, $b=0$, and
$a_k(t)=\alpha_{k+1}t^\gamma$.

\begin{corollary}[Power-law limits]
\label{cor:original}
Let $\gamma>0$, $\alpha_1,\ldots,\alpha_r\ge0$, and $c_m\ge0$ satisfy
$$
 c_m/m^\gamma\to\alpha_{k+1}
 \quad(m\to\infty,\ m\equiv k\pmod r),\quad 0\le k<r.
$$
Let $P_j(z)=z^j$ for $0\le j\le r$,
$P_{n+1}=zP_n-c_{n-r}P_{n-r}$ for $n\ge r$, and
$$
 Q_{n,N}(z)=N^{-n\gamma/d}P_n(N^{\gamma/d}z).
$$
As $N\to\infty$, if $n/N\to t\ge0$, then, locally uniformly off $\Sigma$,
\begin{equation}\label{eq:original-ratio}
 \begin{aligned}
 \frac{Q_{n,N}(z)}{Q_{n+r,N}(z)}&\to\Phi(z,t),\\[7pt]
 1&=\Phi(z,t)\prod_{j=1}^r(z-\alpha_jt^\gamma\Phi(z,t)),
 \quad z^r\Phi(z,t)\to1.
 \end{aligned}
\end{equation}
The normalized zero measures of $Q_{N,N}$ converge to a probability
measure $\mu$ with Cauchy transform
\begin{equation}\label{eq:original-G}
 G(z)=-\frac1r\int_0^1\frac{\partial_z\Phi(z,t)}{\Phi(z,t)}\dd t.
\end{equation}
Put
$$
 A(u)=\prod_{j=1}^r(1-\alpha_j u),
 \quad 1=yF(y)A(F(y)),\quad yF(y)\to1.
$$
For $t>0$, $\Phi(z,t)=zt^{-\gamma}F(z^dt^{-\gamma})$;
$\Phi(z,0)=z^{-r}$. With $y=z^d$,
\begin{equation}\label{eq:original-explicit}
 G(z)=-\frac1{rz}\left(1-\frac d\gamma
 \frac{\displaystyle\int_0^{F(y)}(uA(u))^{1/\gamma}\,\frac{du}{u}}
      {(F(y)A(F(y)))^{1/\gamma}}\right).
\end{equation}
The right-hand side, taken as a whole, is continued from large positive
$z$; this continuation is single-valued because it equals $G$.
For positive $\alpha_j$ and $c_m$, these are the ratio limit and
both transform formulas in
\cite[Theorems~1.1--1.2]{LoureiroVanAssche2026}, for the canonical
initial conditions above.
\end{corollary}

\begin{proof}
The scaled coefficients $c_{m,N}=c_m/N^\gamma$ satisfy
\eqref{eq:bound} and have limiting functions $a_k(t)=\alpha_{k+1}t^\gamma$.
Apply Theorem~\ref{thm:ratios} and Proposition~\ref{prop:compact-zeros}
with $p=r$ and $b=0$. Each product in
\eqref{eq:system} then contains all $r$ phases, giving
\eqref{eq:original-ratio} for $t>0$. If $n/N\to0$, the power
bound $c_m\le K(1+m^\gamma)$ makes all scaled coefficients
up to $n+r$ tend uniformly to zero. The zero bound and
\eqref{eq:ratio-measure} give $\Phi(z,0)=z^{-r}$.
Substitution gives the equation for $F$.
Using the continuation outside the compact support, for sufficiently
large positive $z$ set $u=F(yt^{-\gamma})$ in
\eqref{eq:original-G}; since $t=(yuA(u))^{1/\gamma}$,
$$
 \int_0^1\frac{yt^{-\gamma}F'(yt^{-\gamma})}
                  {F(yt^{-\gamma})}\dd t
 =-\frac{y^{1/\gamma}}\gamma
     \int_0^{F(y)}(uA(u))^{1/\gamma}\,\frac{du}{u}.
$$
This proves \eqref{eq:original-explicit} near infinity.
The full expression agrees there with the Cauchy transform $G$,
so its continuation off the support is single-valued.
\end{proof}

The Appell families in
\cite[Section~4.1, Eq.~(49)]{LoureiroYakubovich2013}
fit the corollary with
$c_m=\binom{m+r}{r}/(r+1)$, $\gamma=r$, and
$\alpha_j=1/((r+1)r!)$.
This coefficient sequence already appears in
\cite[Section~3.1.3]{BenRomdhane2008}.

For $r=2$, $p=4$, and $b=0$, set
$$
 a(t)=\bigl(1+t,\,2+t^2,\,3+\sqrt t,\,1+t^{-1}\bigr),
 \qquad c_{4m+k,N}=a_k\bigl((4m+k+1)/N\bigr),\quad 0\le k<4.
$$
All coefficients are positive, and both \eqref{eq:local} and
\eqref{eq:measure} hold. Here $p\ne r$, the profiles have no common
power form, and $a_3$ is not integrable at zero. Moreover,
$c_{3,N}=1+N/4$, so \eqref{eq:bound} fails. By
\eqref{eq:periodic-moments},
$\int\zeta^3\dd\mu_{a(t),0}(\zeta)=\frac34\sum_{k=0}^3a_k(t)$.
Since $\zeta^3\ge0$ on the three-ray star, Tonelli's theorem gives
$\int\zeta^3\dd\mu_T(\zeta)=\infty$; thus $\mu_T$ has noncompact support.

For $r=p=1$, take $b=0$, $c_{m,N}=1$ except $c_{N-1,N}=2$.
Condition~\eqref{eq:L1} holds with $a_0=1$, $h_0=0$.
Nevertheless, writing $R_{n,N}=q_{n,N}/q_{n+1,N}$ and
$f(z)=(z-\sqrt{z^2-4})/2$, with $\sqrt{z^2-4}\sim z$, gives
$$
 R_{N-1,N}(z)\to f(z),\quad
 R_{N,N}(z)=\frac1{z-2R_{N-1,N}(z)}\to\frac1{\sqrt{z^2-4}}.
$$
The first ratio belongs to the constant recurrence. The two limits
are different although both scaled indices tend to $1$.
Thus the zero theorem applies even when the ratios have no common limit
along all sequences with $n/N\to1$.

The zero theorem also allows functions that are not integrable: for $p=1$, $b=0$ and
$c_{m,N}=N/(m+1)$ the limiting function is $a(t)=1/t$, which is not in
$L^1(0,T)$. Convergence in measure holds, and the limiting
measure has infinite $d$th moment.
Weak $L^1$ convergence alone, however, is insufficient.
For $r=p=1$ and $b=0$, the step functions for coefficients alternating
between $1$ and $3$ converge weakly in $L^1$ to $a(t)=2$, but
their limiting fourth zero moment is $22$; the constant recurrence
with $c=2$ has fourth moment $24$.

For $p=d$, $b=0$, and only the zeroth coefficient phase nonzero,
$$
 q_{dm+s,N}(z)=z^s\prod_{j=0}^{m-1}(z^d-c_{dj,N}),\quad m\ge0,\quad 0\le s<d.
$$
Hence $\mu_{(u,0,\ldots,0),0}=\rho_u$, where $\rho_u$ is uniform
on the roots of $z^d=u$ for $u>0$ and $\rho_0=\delta_0$.
Given a probability $\nu$ on $[0,\infty)$, take $a_0$ to be its
generalized quantile on $(0,1)$ and zero for $t=0$ and $t\ge1$, and set
$c_{dm,N}=a_0(dm/N)$, with all other phases zero.
The step functions converge almost everywhere, and thus in measure.
At $T=1$ this gives $\mu_1=\int\rho_u\dd\nu(u)$.
Every rotation-invariant probability on $\Sigma$ can therefore occur.

\section{Proofs}

We first note that the joint limits in \eqref{eq:local} force continuity.
Fix a phase $k$ and a sequence $t_j\to t>0$. For each $j$, choose
$N_j\to\infty$ and $m_j\equiv k\pmod p$ such that
$|m_j/N_j-t_j|<1/j$ and
$|c_{m_j,N_j}-a_k(t_j)|+|b_{m_j,N_j}-h_k(t_j)|<1/j$.
This is possible by the joint limit at each fixed $t_j$. Since
$m_j/N_j\to t$, \eqref{eq:local} gives $a_k(t_j)\to a_k(t)$ and
$h_k(t_j)\to h_k(t)$.

We now keep $N$ fixed and omit it from the notation.
Known interlacing results
\cite[Theorem~2.7 and Section~4.2]{DelvauxLopez2013} and partial fractions
\cite[proof of Lemma~2.2]{KuijlaarsVanAssche1999} give
\begin{equation}\label{eq:ratio-measure}
 R_n(z)=\frac{q_n(z)}{q_{n+1}(z)}
       =\int\frac{\dd\lambda_n(\zeta)}{z-\zeta},
\end{equation}
where $\lambda_n$ is a probability on the zeros of $q_{n+1}$ in
$\Sigma$. For $r\ge2$, write $q_n=z^jP(u)$, where $u=z^d$
and $j=n\bmod d$. Then $R_n=z^rP(u)/(uQ(u))$ if $j<r$,
and $R_n=z^rP(u)/Q(u)$ if $j=r$. Interlacing gives nonnegative
residues; they sum to one since each reduced ratio is
$u^{-1}+O(u^{-2})$ at infinity.
For $r\ge2$, the identity
$$
 \frac1d\sum_{j=0}^{d-1}\frac1{z-\omega^ju^{1/d}}
 =\frac{z^r}{z^d-u}
$$
gives a rotation-invariant probability measure
$\lambda_n$; the zero measures are also rotation invariant.
The representation for nonnegative coefficients follows by taking
$c_m+\varepsilon$ and letting $\varepsilon\downarrow0$.

The zeros of $q_{n+1}$ are the eigenvalues of its recurrence matrix,
with diagonal $b_m$, superdiagonal $1$ and $r$th subdiagonal $c_m$.
If $|b_m|,c_m\le C$ for $0\le m\le n$, conjugation by
$\operatorname{diag}(1,s,\ldots,s^n)$ bounds the maximum absolute row sum
by $s+Cs^{-r}$, with an additional $C$ when $r=1$.
Minimizing over $s>0$ gives \eqref{eq:root-bound}; $C=0$ follows
by continuity. Thus the zeros lie in $\Sigma_{B_r(C)}$.
The bound on $|R_n(z)|$ is then
$\operatorname{dist}(z,\Sigma_{B_r(C)})^{-1}$.

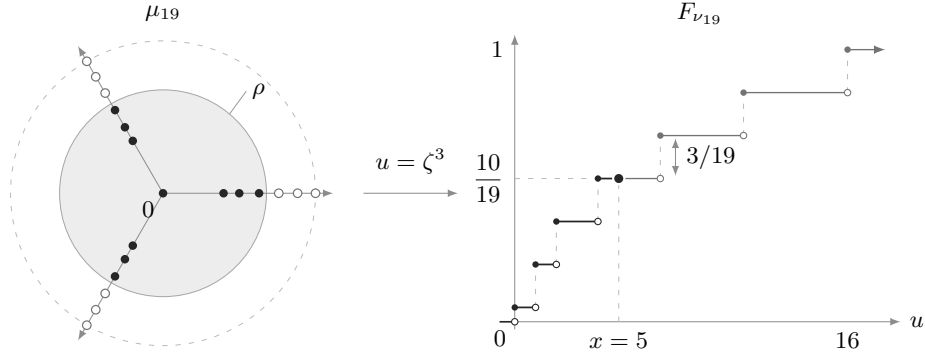
\begin{figure}[!htb]
\centering
\begin{tikzpicture}[font=\small,>=latex,x=1cm,y=1cm,
  math axis/.style={->,black!45,line width=0.35pt},
  math guide/.style={black!27,line width=0.3pt,
                    dash pattern=on 2pt off 3pt},
  count curve/.style={black!55,line width=0.5pt},
  marked curve/.style={black!85,line width=0.6pt},
  counted zero/.style={fill=black!85},
  other zero/.style={draw=black!55,fill=white,line width=0.45pt}]
  \pgfmathsetmacro{\cutRadius}{0.8*pow(5,1/3)}
  \pgfmathsetmacro{\outerRadius}{0.8*pow(16,1/3)}
  \node at (0,2.38) {$\mu_{19}$};
  \draw[math guide] (0,0) circle (\outerRadius);
  \fill[black!7] (0,0) circle (\cutRadius);
  \draw[black!35,line width=0.35pt] (0,0) circle (\cutRadius);
  \draw[black!40,line width=0.3pt]
    (50:\cutRadius) -- (1.1,1.35);
  \node[right,inner sep=2pt] at (1.1,1.35) {$\rho$};
  \foreach \angle in {0,120,240}{
    \draw[math axis] (0,0) -- (\angle:2.25);
    \foreach \u in {1,2,4}{
      \pgfmathsetmacro{\radius}{0.8*pow(\u,1/3)}
      \path[counted zero] (\angle:\radius) circle (1.6pt);
    }
    \foreach \u in {7,11,16}{
      \pgfmathsetmacro{\radius}{0.8*pow(\u,1/3)}
      \draw[other zero] (\angle:\radius) circle (1.6pt);
    }
  }
  \path[counted zero] (0,0) circle (1.6pt);
  \node[below left,inner sep=3pt] at (0,0) {$0$};
  \draw[math axis] (2.65,0) -- (3.9,0)
    node[midway,above,yshift=4pt,text=black] {$u=\zeta^3$};
  \begin{scope}[shift={(4.65,-1.7)}]
    \node at (2.45,4.08) {$F_{\nu_{19}}$};
    \draw[math axis] (-0.2,0) -- (5.1,0)
      node[right,text=black] {$u$};
    \draw[math axis] (0,-0.1) -- (0,3.8);
    \draw[math guide] (0,{3.6*10/19}) -- (1.375,{3.6*10/19});
    \draw[math guide] (1.375,0) -- (1.375,{3.6*10/19});
    \node[left,inner sep=4pt] at (0,{3.6*10/19})
      {$\displaystyle\frac{10}{19}$};
    \node[left,inner sep=4pt] at (0,3.6) {$1$};
    \node[below left,inner sep=3pt] at (0,0) {$0$};
    \node[below,inner sep=4pt] at (1.375,0) {$x=5$};
    \node[below,inner sep=4pt] at (4.4,0) {$16$};
    \draw[marked curve] (-0.2,0) -- (0,0);
    \foreach \u/\count/\nextU in {0/1/1,1/4/2,2/7/4,4/10/7,7/13/11,11/16/16}{
      \draw[count curve] ({0.275*\u},{3.6*\count/19})
        -- ({0.275*\nextU},{3.6*\count/19});
      \draw[math guide] ({0.275*\nextU},{3.6*\count/19})
        -- ({0.275*\nextU},{3.6*(\count+3)/19});
    }
    \draw[count curve,->] (4.4,3.6) -- (4.9,3.6);
    \foreach \u/\count/\nextU in {0/1/1,1/4/2,2/7/4,4/10/5}{
      \draw[marked curve] ({0.275*\u},{3.6*\count/19})
        -- ({0.275*\nextU},{3.6*\count/19});
    }
    \foreach \u/\count in {0/0,1/1,2/4,4/7}{
      \draw[draw=black!85,fill=white,line width=0.4pt]
        ({0.275*\u},{3.6*\count/19}) circle (1.2pt);
    }
    \foreach \u/\count in {0/1,1/4,2/7,4/10}{
      \fill[black!85] ({0.275*\u},{3.6*\count/19}) circle (1.2pt);
    }
    \foreach \u/\count in {7/10,11/13,16/16}{
      \draw[draw=black!55,fill=white,line width=0.4pt]
        ({0.275*\u},{3.6*\count/19}) circle (1.2pt);
      \fill[black!55] ({0.275*\u},{3.6*(\count+3)/19}) circle (1.2pt);
    }
    \fill[white] (1.375,{3.6*10/19}) circle (2.5pt);
    \fill[black!85] (1.375,{3.6*10/19}) circle (1.7pt);
    \draw[<->,black!45,line width=0.3pt]
      (2.125,{3.6*10/19+0.035}) -- (2.125,{3.6*13/19-0.035})
      node[midway,right,inner sep=4pt,text=black] {$3/19$};
  \end{scope}
\end{tikzpicture}
\caption{For $r=2$, set $c_{3j}=1+j(j+1)/2$ for $0\le j<6$,
and all other coefficients to zero. Each zero of $q_{19}$ has mass $1/19$.
The real image $\nu_{19}$ has $F_{\nu_{19}}(5)=10/19$, the mass of the ten
zeros in the shaded disk of radius $\rho=5^{1/3}$.
Each nonzero orbit adds~$3/19$. The outer circle is a guide.}
\label{fig:star-image}
\end{figure}

Moments determine probabilities on a compact real interval.
A rotation-invariant probability measure on $\Sigma_B$ is determined by its
image on $[0,B^d]$ under $\zeta\mapsto\zeta^d$, since its
weights on each orbit are equal.
Figure~\ref{fig:star-image} shows the same normalized count on the star
and on the real axis; each nonzero orbit in the figure contributes
the weight $d/n$ used in the next lemma.

\begin{samepage}
\begin{lemma}[Changing recurrence coefficients]\label{lem:counting}
Let $n\ge1$ and compare two recurrences satisfying
\eqref{eq:recurrence} and the stated sign conditions. Set
$$
 K=\#\{0\le j<n:b_j\ne\widetilde b_j\}
   +\#\{0\le j<n-r:c_j\ne\widetilde c_j\}.
$$
Let $\mu_n,\widetilde\mu_n$ be their normalized zero measures, and put
$\nu_n=\mu_n$, $\widetilde\nu_n=\widetilde\mu_n$ when $r=1$;
when $r\ge2$, take their images under $\zeta\mapsto\zeta^d$.
For a real probability measure $\nu$, write
$F_\nu(x)=\nu((-\infty,x])$. Then
\begin{equation}\label{eq:counting}
 \sup_{x\in\R}|F_{\nu_n}(x)-F_{\widetilde\nu_n}(x)|
 \le\begin{cases}K/n,&r=1,\\[7pt] dK/n,&r\ge2.\end{cases}
\end{equation}
\end{lemma}
\end{samepage}

\begin{proof}
First take all $c_j>0$. For $r=1$ the zeros are real
and simple. For $r\ge2$, write $q_n(z)=z^sP(z^d)$, $0\le s<d$;
the zeros of $P$ are positive and simple by
\cite[Theorem~2.2]{BenRomdhane2008}; see also
\cite[Section~4.2]{DelvauxLopez2013}.
Interpolate one coefficient between its two endpoint values.
The subdiagonal coefficients stay positive throughout this interval, so
the relevant real roots remain simple and can be followed continuously.
The recurrence determinant is affine in each coefficient, which
occupies one matrix entry. For fixed real $x$, the value $q_n(x)$
when $r=1$, or $P(x)$ otherwise, is affine in the parameter and has
at most one zero, unless it vanishes identically. Since the zeros are
simple, at most one zero can pass through $x$ as the coefficient varies.
If this value vanishes identically as a function of the parameter,
a zero remains at $x$ and no other zero can reach it.
Each positive root of $P$ has weight $d/n$ in the real image.
Successive changes give the stated bounds. For nonnegative coefficients,
replace all $c_j$ by $c_j+\varepsilon$ and pass to the limit at
continuity points, then use right continuity.
\end{proof}

Consequently, changing $o(n)$ recurrence coefficients by arbitrary
amounts, while preserving the sign conditions, does not alter any
weak limit of the normalized zero measures.

\subsection{A calculation with power series}
Set $v=z^{-1}$ and $S_{n,N}=zq_{n,N}/q_{n+1,N}$,
with $b_{m,N}=c_{m,N}=0$ and $S_{m,N}=1$ for $m<0$. Then
\begin{equation}\label{eq:S}
 S_{n,N}(v)=\left(1-b_{n,N}v-
 c_{n-r,N}v^d\prod_{j=1}^rS_{n-j,N}(v)\right)^{-1}.
\end{equation}
Induction on $\ell$ shows that $[v^\ell]S_{n,N}$ and
$[v^\ell]\log S_{n,N}$ are polynomials in recurrence coefficients
with indices in $[n-r\ell,n]$: preceding series enter only at
lower orders. The form of these polynomials does not depend on $n$ or $N$,
and they are Lipschitz on bounded sets.

For fixed vectors $a,h$, define formal power series $f_k(v;a,h)$
with constant term $1$ by
\begin{equation}\label{eq:f}
 f_k(v;a,h)=\left(1-h_kv-a_{k-r}v^d
                  \prod_{j=1}^rf_{k-j}(v;a,h)\right)^{-1},
\end{equation}
and put $\Psi(v;a,h)=\prod_{k=0}^{p-1}f_k(v;a,h)$.
Recursion by degree determines their coefficients uniquely.
For the fixed periodic recurrence,
$[v^\ell]S_n=[v^\ell]f_k$ once $n\ge r\ell$ and $n\equiv k\pmod p$.
For $C=\max_k\{a_k,|h_k|\}$, the root bound and
\eqref{eq:ratio-measure} give $|[v^\ell]f_k|\le B_r(C)^\ell$.
Thus these series define the unique analytic branch near zero.

\subsection{Proof of the ratio theorem}
Fix $n/N\to t>0$ with $n\equiv k\pmod p$.
Each coefficient at a fixed distance from $n$ converges to the value
of its limiting function at $t$. The coefficient calculation above gives
$[v^\ell]S_{n,N}\to[v^\ell]f_k(v;a(t),h(t))$ for every $\ell$.
These are the moments of the probability measures in
\eqref{eq:ratio-measure}, now denoted by $\lambda_{n,N}$.
For fixed periodic coefficients, the moments eventually remain
unchanged in each phase. Since the measures have a common compact
support, they have a unique weak limit in each phase.
For the varying coefficients, the second moment of $\lambda_{n,N}$
is $b_{n,N}^2+c_{n-1,N}$ when $r=1$.
When $r\ge2$, its image under $\zeta\mapsto\zeta^d$ is a probability
measure on $[0,\infty)$ with first moment $c_{n-r,N}$.
These bounds give tightness; higher even moments, or higher moments
after applying $\zeta\mapsto\zeta^d$, give uniform integrability.
Indeed, for $r=1$, choose an integer $q$ with $2q>\ell$. Then
$$
 \int_{|x|>R}|x|^\ell\dd\lambda_{n,N}(x)
 \le R^{\ell-2q}\int|x|^{2q}\dd\lambda_{n,N}(x).
$$
For $r\ge2$, writing $\eta_{n,N}$ for the real image, one has
$$
 \int_{u>R}u^m\dd\eta_{n,N}(u)
 \le R^{-1}\int u^{m+1}\dd\eta_{n,N}(u).
$$
The moments on the right are uniformly bounded along the sequence by
the coefficient calculation.
Every weak limit has the same moments as the limiting probability
measure in phase $k$ for the periodic coefficients $a(t),h(t)$.
The geometric bound on these moments forces compact support: if
$B=B_r(\max_j\{a_j(t),|h_j(t)|\})$, then any weak limit $\lambda$
in the real case satisfies, for $R>B$,
$\lambda(\{|x|\ge R\})\le(B/R)^{2m}$ for every integer $m\ge1$;
any weak limit $\eta$ of the image measures satisfies
$\eta([R,\infty))\le(B^d/R)^m$ when $R>B^d$.
Letting $m\to\infty$ proves the assertion.
Moment uniqueness then identifies the weak limit $\lambda_{k,t}$.
The ratios converge to
$\varphi_k(z,t)=\int(z-\zeta)^{-1}\dd\lambda_{k,t}(\zeta)$
locally uniformly off $\Sigma$ by boundedness and equicontinuity
of the functions $(z-\zeta)^{-1}$.
Near infinity, $\varphi_k=z^{-1}f_k(z^{-1};a(t),h(t))$.

For $r=1$, the imaginary part of $\varphi_k$ has a strict sign off $\R$.
For $r\ge2$, rotation invariance gives
$\varphi_k=z^r\int(z^d-u)^{-1}\dd\eta_{k,t}(u)$ for a probability
$\eta_{k,t}$, compactly supported on $[0,\infty)$.
The imaginary part of the integral is nonzero if $z^d$ is nonreal;
the integral is negative if $z^d<0$.
Thus the limits are nonzero, and passing to the limit in the
recurrence gives \eqref{eq:system}.
Products of $p$ consecutive ratios give the limit of $q_{n,N}/q_{n+p,N}$.
Uniqueness of the formal power series and the identity theorem determine
the analytic branch on
each component off $\Sigma$.

The measures $\lambda_{k,t}$ depend weakly continuously on the
coefficient vectors. Their moments are polynomials in the entries of
these vectors, and their supports lie in a common compact set when
the entries remain bounded. If uniform convergence on a compact
subinterval of $(0,\infty)$ failed, we could select indices witnessing
failure, extract a fixed phase and a subsequence with $n/N\to t$,
and contradict the sequential limit together with weak continuity of
the periodic measures. Cauchy's estimates then give convergence of
every fixed derivative.
\qed

\begin{proposition}[Local ratio limits]\label{prop:local-ratios}
Under \eqref{eq:recurrence} and the stated sign conditions, let
$n_N\to\infty$, $a\in[0,\infty)^p$ and $h\in\R^p$,
with $h=0$ if $r\ge2$. The phases here are relative to $n_N$;
no congruence condition on $n_N$ is imposed. Suppose that
$$
 \begin{aligned}
 c_{n_N+j,N}&\to a_{j\bmod p}\quad\text{for every fixed integer }j\le p-r-1,\\[7pt]
 b_{n_N+j,N}&\to h_{j\bmod p}\quad\text{for every fixed integer }j\le p-1.
 \end{aligned}
$$
Then, locally uniformly on $\C\setminus\Sigma$,
$$
 \frac{q_{n_N,N}(z)}{q_{n_N+1,N}(z)}\to\varphi_0(z;a,h),
 \quad
 \frac{q_{n_N,N}(z)}{q_{n_N+p,N}(z)}\to\Phi(z;a,h),
$$
also after any fixed number of derivatives in $z$.
\end{proposition}

\begin{proof}
The hypotheses control each fixed finite backward window, without
imposing a global coefficient bound. For $0\le j<p$, the coefficient
calculation gives
$[v^\ell]S_{n_N+j,N}\to[v^\ell]f_{j\bmod p}(v;a,h)$;
the same moment argument proves convergence of these shifted ratios,
and multiplication gives the second limit. The highest coefficient indices
entering $q_{n_N+p,N}$ are $n_N+p-r-1$ for $c$ and $n_N+p-1$ for $b$.
\end{proof}

\subsection{Proof of Proposition~\ref{prop:compact-zeros}}
Put $L_\ell(n,N)=[v^\ell]\log S_{n,N}$.
The identity $q_{n,N}=z^n\prod_{j<n}S_{j,N}(z^{-1})^{-1}$ gives
\begin{equation}\label{eq:moment}
 \int\zeta^\ell\dd\mu_{n,N}(\zeta)
 =\frac\ell n\sum_{j<n}L_\ell(j,N),\quad\ell\ge1.
\end{equation}
For fixed periodic coefficients, $L_\ell(n)$ eventually remains
unchanged in each phase. Averaging in \eqref{eq:moment} and using
the uniqueness from moments proved above establishes the existence
of $\mu_{a,h}$, with
\begin{equation}\label{eq:periodic-moments}
 \int\zeta^\ell\dd\mu_{a,h}(\zeta)
 =\frac\ell p[v^\ell]\log\Psi(v;a,h).
\end{equation}
The moments are polynomials in $a,h$, and the supports lie in a
common compact set when $a,h$ remain bounded. This proves weak
continuity in $a,h$, so the average in \eqref{eq:zero-average} is well defined.

Fix $T<U<L$. The step coefficients on $(0,U)$ are eventually
bounded by $C_L$, as are their limits almost everywhere; choose
bounded representatives. The zero measures of $q_{n,N}$ and the
measures $\mu_{a(t),h(t)}$ are supported on $\Sigma_B$, with $B=B_r(C_L)$.
Extend the step functions and their limits by zero below $0$.
Each recurrence coefficient entering
$L_\ell(p\lfloor Nt/p\rfloor+k,N)$ comes from one of the step functions
$A_{j,N}$ or $H_{j,N}$, translated by a fixed multiple of $p/N$.
Continuity of translations in $L^1$ and the Lipschitz estimates
for the coefficient polynomials give
$$
 L_\ell(p\lfloor Nt/p\rfloor+k,N)
 \to[v^\ell]\log f_k(v;a(t),h(t))
 \quad\text{in }L^1(0,U).
$$
Each complete interval of length $p/N$ corresponds to $p$ consecutive summands.
The incomplete interval and the replacement of $n/N$ by $T$ contribute
an error of $O(N^{-1}+|n/N-T|)$, since the summands are bounded. Thus \eqref{eq:moment} gives
\begin{equation}\label{eq:average-moments}
 \int\zeta^\ell\dd\mu_{n,N}(\zeta)\to
 \frac\ell{pT}\int_0^T[v^\ell]\log\Psi(v;a(t),h(t))\dd t.
\end{equation}
These are the moments of $\mu_T$ by \eqref{eq:periodic-moments}.
The uniqueness from moments proved above gives \eqref{eq:zero-average}.

For $\max_k\{a_k,|h_k|\}\le C_L$, define on $\C\setminus\Sigma_B$
\begin{equation}\label{eq:potential}
 \Phi(z;a,h)=z^{-p}\exp\left\{-p\int
       \operatorname{Log}(1-\zeta/z)\dd\mu_{a,h}(\zeta)\right\}.
\end{equation}
The principal logarithm is analytic here: its cut would require
$z\in[0,\zeta]\subset\Sigma_B$.
By \eqref{eq:periodic-moments}, this nonzero function equals
$z^{-p}\Psi(z^{-1};a,h)$ near infinity, so it analytically continues
the product $\prod_{k=0}^{p-1}\varphi_k$.
Differentiation gives
$$
 -\frac1p\frac{\partial_z\Phi}{\Phi}
 =\int\frac{\dd\mu_{a,h}(\zeta)}{z-\zeta}.
$$
Averaging proves \eqref{eq:G}. The common compact support justifies
differentiation and Fubini; weak convergence gives locally uniform
convergence of the Cauchy transforms by the same boundedness
and equicontinuity argument used for the ratios.
\qed

\section{The power-law case: estimates and earlier laws}

We now return to the power-law specialization in
Corollary~\ref{cor:original}. The following observations concern
estimates, measure identification, and earlier explicit laws in
\cite{LoureiroVanAssche2026}.

\paragraph{A ratio bound:}
For $n\ge r$ and $z\in\C\setminus\Sigma$, a valid replacement for the first bound in
\cite[(2.10)]{LoureiroVanAssche2026} is
$$
 \left|\frac{Q_{n-r,N}(z)}{Q_{n+1,N}(z)}\right|
 \le\delta^{-(r+1)},\quad \delta=\operatorname{dist}(z,\Sigma),
$$
obtained by multiplying $r+1$ bounds from \eqref{eq:ratio-measure}.
On each fixed compact set off the star, this changes only the
constants used in the convergence estimates.

\paragraph{The star measure and its real image:}
In the symmetric case, the probability
$\nu=(\zeta\mapsto\zeta^d)_*\mu$ satisfies
$$
 G_\mu(z)=z^rG_\nu(z^d),\quad
 G_\nu(y)=\int\frac{\dd\nu(u)}{y-u}.
$$
This standard identity appears in
\cite[Proposition~2.2, $k=j$]{LopezGarciaMinaDiaz2018}.
The function $\widehat G$ on p.~19 of the preprint instead has
mass $1/d$ and must be multiplied by $d$ to represent $\nu$.
With this convention the density, where defined, is
$((G_\nu)_--(G_\nu)_+)/(2\pi i)$.
The real-line inversion applies to $\nu$ and has this jump order;
the Cauchy integral itself has no extra factor $1/(2\pi)$.

\paragraph{The ordinary $S$-transform:}
For the positive parameters in the preprint, formula~(3.8) concerns
$\nu$, not the original star measure $\mu$.
Rotational symmetry gives $\int\zeta\dd\mu=0$, so the local
analytic inverse in the ordinary $S$-transform is unavailable
for $\mu$; see \cite[Section~4.5.3]{MingoSpeicher2017} for the convention.
For $r=1$, $c_j=j+1$, the limiting star measure is semicircular
on $[-2,2]$. Its image under $x\mapsto x^2$ is Marchenko--Pastur
on $[0,4]$, with $S(w)=1/(1+w)$, the expression in (3.8).
The probability measures on $[0,\infty)$ in that discussion therefore
describe the powers of the zeros.

\paragraph{Earlier limiting laws:}
For the case $\gamma=1$ in \cite[Section~4.3]{LoureiroVanAssche2026},
with $r\ge2$, $\alpha_1=\cdots=\alpha_{r-1}=\alpha>0$ and
$\alpha_r=(r+1)\alpha$, the real-image law $\nu$ is the dilation
by $\alpha((r+1)/r)^r$ of $J_{r,r-1,r}/r$, where $J_{r,s,a}$
is the measure of mass $a$ studied in
\cite[Theorem~2.3 and Proposition~2.4]{GawronskiNeuschelStivigny2016}.
The branch points of the corresponding algebraic curve are special
cases of the formulas there.
The same law occurs as the step-line zero limit in
\cite[Theorem~3.14 and Corollary~3.15]{Wolfs2024}:
take $p=r$, $q=r-1$ there and dilate by
$\alpha(r+1)^r/r^{r-1}$.
For $r=2$, it is the dilation by $9\alpha/2$ of the equal-proportion
law in \cite{VanAsscheWolfs2023}, whose density is explicit.
The polar parametrization follows from the same quadratic
calculation used in \cite{Neuschel2020} for $J_{r,r-1,1}$,
applied to $J_{r,r-1,r}$. Stieltjes inversion gives the density,
and local expansions give the endpoint constants.
These formulas follow from the known law.
Theorems~\ref{thm:ratios} and~\ref{thm:zeros} extend the coefficient
hypotheses under which convergence holds.

We likewise offer our own arguments for scrutiny by readers and by
future, more capable large language models, and welcome any corrections.

\section{Development of the results and use of AI}
\label{sec:ai-use}

ChatGPT 5.6 Sol drafted the proofs of the ratio theorem and the bounded-zero result separately from photographs of the authors’ notes. It identified relevant earlier works and helped to organize ideas dispersed throughout the notes. ChatGPT 6 assisted with subsequent reviews, in which it identified no errors, proposed the local assumption for the ratios, and drafted Appendix~\ref{sec:unbounded}, which is entirely its own contribution.

We address these results because they appeared in a September 2026 preprint,
allowing us to contribute to the discussion before publication. AI also suggested extensions and further results. We see little value in adding them to the
literature merely because they can now be developed more quickly, unless
they address a structural question or a concrete mathematical or social need.

As authors, we take full responsibility for the entire content of this work.

\section*{Acknowledgements}

The first author
acknowledges financial support from the Centre for Mathematics of
the University of Coimbra (CMUC), funded by the Portuguese Foundation for
Science and Technology (FCT), under the projects UID/00324/2025
(\url{https://doi.org/10.54499/UID/00324/2025}) and UID/PRR/00324/2025. The
first author also acknowledges financial support from the FCT under the grant
\url{https://doi.org/10.54499/2022.00143.CEECIND/CP1714/CT0002}.

\appendix
\section{Proof of Theorem~\ref{thm:zeros}}
\label{sec:unbounded}

We remove the coefficient bounds by truncation, using
Proposition~\ref{prop:compact-zeros} and Lemma~\ref{lem:counting}.

\begin{proof}
For $M>0$, truncate the coefficients by
$$
 c^M=\min\{c,M\},\quad
 b^M=\max\{-M,\min\{b,M\}\},
$$
and denote the corresponding truncations of $a,h$ by $a^M,h^M$.
The truncated step functions converge in measure and are uniformly
bounded, so they converge in $L^1(0,U)$ for every $U>0$.
Proposition~\ref{prop:compact-zeros} applies to the recurrences with truncated coefficients.
If $K_N(M)$ is the sum of the numbers of changed diagonal and
subdiagonal coefficients entering $q_{n,N}$, use
$$
 \{A_{k,N}>M\}\subset
 \{a_k>M/2\}\cup\{|A_{k,N}-a_k|>M/2\},
$$
and the analogous inclusion for $|H_{k,N}|$ and $|h_k|$.
Each phase coefficient accounts for an interval of length $p/N$;
the initial and final incomplete intervals contribute $O(1)$ coefficients.
Replacing $n/N$ by $T$ changes the normalized count by
$O(|n/N-T|+N^{-1})$. Convergence in measure therefore gives
$$
 \limsup_{N\to\infty}\frac{K_N(M)}n
 \le\frac1{pT}\sum_{k=0}^{p-1}
 \bigl(|\{a_k>M/2\}\cap(0,T)|
       +|\{|h_k|>M/2\}\cap(0,T)|\bigr).
$$
The right-hand side tends to zero as $M\to\infty$ because the
functions $a_k,h_k$ are finite almost everywhere.
Lemma~\ref{lem:counting} therefore makes the distance between the
original and truncated real measures arbitrarily small in the
L\'evy metric, using their real images when $r\ge2$.
For almost every $t$, sufficiently large $M$ leaves every component
of $a(t)$ and $h(t)$ unchanged. Dominated convergence for bounded
continuous test functions shows that the averages
$\frac1T\int_0^T\mu_{a^M(t),h^M(t)}\dd t$ converge weakly to $\mu_T$.
The triangle inequality in the L\'evy metric gives \eqref{eq:zero-average}.
For $r\ge2$, weak convergence of the real images gives weak
convergence on the star by averaging each test function over the
$d$ roots of $u\ge0$.

Finally, the functions $z\mapsto(z-\zeta)^{-1}$ are uniformly bounded
and equicontinuous on compact subsets of $\C\setminus\Sigma$,
uniformly for $\zeta\in\Sigma$. Thus weak convergence gives the stated
locally uniform convergence. For every pair of coefficient vectors $a,h$
with finite entries,
\eqref{eq:potential} gives
$-\partial_z\Phi/(p\Phi)=G_{\mu_{a,h}}$.
The bound $|G_{\mu_{a,h}}(z)|\le\operatorname{dist}(z,\Sigma)^{-1}$
justifies averaging this identity without any integrability
assumption on the functions $a_k,h_k$.
\end{proof}

\section{The power-law case from earlier identities}
\label{sec:power-known}

Assume the positive coefficients and parameters of
Corollary~\ref{cor:original}, and set $q_{n,N}=Q_{n,N}$.
The power-law setting of the two limits in
\cite[Theorems~1.1--1.2]{LoureiroVanAssche2026} is more elementary
than the general coefficient regimes considered here: the common
scaling form and the bounds on finite scaled intervals permit a
particularly simple proof from earlier continued-fraction identities.
Interlacing, normality, and moment convergence justify the passage
from formal identities to analytic and weak limits.
Let $H$ be the unique analytic function near zero, with $H(0)=1$,
given by the known period-$r$ identity
\cite[Eq.~(10.34)]{PetreolleSokalZhu2023}
$$
 H(x)=\prod_{j=1}^r(1-\alpha_jxH(x))^{-1}.
$$

To read the finite branched-fraction identity
\cite[Eq.~(2.27)]{PetreolleSokalZhu2023} from the last recurrence
coefficient, fix $n$, put $x=v^d$, and set
$g_k(x)=S_{n-k,N}(v)$ for $k\ge0$. With
$\beta_{k+r}=c_{n-r-k}/N^\gamma$ and $c_j=0$ for $j<0$,
\eqref{eq:S} becomes
$$
 g_k(x)=\left(1-\beta_{k+r}x\prod_{j=1}^rg_{k+j}(x)\right)^{-1}.
$$
The zero coefficients terminate this continued fraction.
Its coefficient $[v^{dm}]S_{n,N}$ is a homogeneous polynomial
of degree $m$ in finitely many neighbouring coefficients $c_j/N^\gamma$.
Consequently, when $n/N\to t>0$, the Laurent coefficients of
$Q_{n,N}/Q_{n+r,N}$ tend to those of
$\Phi(z,t)=z^{-r}H(t^\gamma z^{-d})$.
The bound $c_j\le K(1+j^\gamma)$ gives a common compact star
whenever $n/N$ stays bounded. Positive partial fractions
\eqref{eq:ratio-measure}, obtained from interlacing
\cite[Theorem~2.2 and its proof]{BenRomdhane2008}, give normality
off $\Sigma$; the coefficient limits then give local uniform
convergence there. The limit is nonzero by Hurwitz's theorem.
If $n/N\to0$, all scaled coefficients up to
$n+r$ tend uniformly to zero, giving $\Phi(z,0)=z^{-r}$.

For the zeros, use
$Q_{N,N}=z^N\prod_{j<N}S_{j,N}(z^{-1})^{-1}$.
The coefficients $[v^{dm}]\log S_{j,N}$ are likewise homogeneous
polynomials of degree $m$. On each interval $[\varepsilon,1]$,
$\varepsilon>0$, the coefficient limits are uniform in each phase;
the first $\varepsilon N$ summands contribute $O(\varepsilon)$ after
normalization, by the common coefficient bound. Averaging over the
$r$ residue classes, taking Riemann sums, and then letting
$\varepsilon\downarrow0$ gives
$$
 \begin{aligned}
 \lim_{N\to\infty}\int\zeta^{dm}\dd\mu_{N,N}(\zeta)
 &=\frac{dm}{r(1+\gamma m)}[x^m]\log H(x)\\[7pt]
 &=\frac{d}{r(1+\gamma m)}
   [u^m]\prod_{j=1}^r(1-\alpha_j u)^{-m},\quad m\ge1.
 \end{aligned}
$$
The last equality is Lagrange inversion; moments of positive
integer order not divisible by $d$ vanish by rotational symmetry.
Rotation invariance, the common compact support, and moment uniqueness
for the real image under $\zeta\mapsto\zeta^d$ determine the weak limit.
Its Cauchy transform is \eqref{eq:original-G}, as seen by expanding
$-r^{-1}\int_0^1\partial_z\log\Phi(z,t)\dd t$ at infinity;
the substitution in the corollary gives \eqref{eq:original-explicit}.

This deduction uses bounds on finite scaled intervals and a common
power form. Our main results allow arbitrary periods and independent limiting
functions. For zero limits, convergence in measure suffices,
without coefficient bounds or integrability assumptions.

\makeatletter
\renewcommand{\@biblabel}[1]{\@defaultbiblabelstyle{#1}}
\def\bibsetup{}
\makeatother

\end{document}